# Exchangeable lower previsions

GERT DE COOMAN[1,*], ERIK QUAEGHEBEUR[1,**] and
ENRIQUE MIRANDA[2]

[1]*Ghent University, SYSTeMS Research Group, Technologiepark–Zwijnaarde 914, 9052 Zwijnaarde, Belgium. E-mails:* [*]*gert.decooman@ugent.be;* [**]*erik.quaeghebeur@ugent.be*
[2]*University of Oviedo, Dept. of Statistics and Operations Research. C-Calvo Sotelo, s/n, 33007, Oviedo, Spain. E-mail: mirandaenrique@uniovi.es*

We extend de Finetti's [*Ann. Inst. H. Poincaré* **7** (1937) 1–68] notion of exchangeability to finite and countable sequences of variables, when a subject's beliefs about them are modelled using coherent lower previsions rather than (linear) previsions. We derive representation theorems in both the finite and countable cases, in terms of sampling without and with replacement, respectively.

## 1. Introduction

This paper deals with belief models for both finite and countable sequences of exchangeable random variables taking a finite number of values. When such sequences of random variables are assumed to be exchangeable, this more-or-less means that the specific order in which they are observed is deemed irrelevant.

The first detailed study of exchangeability was made by de Finetti [5] (with the terminology of 'equivalent' events). He proved the now famous representation theorem, which is often interpreted as stating that a sequence of random variables is exchangeable if it is conditionally independent and identically distributed (i.i.d.). Other important work on exchangeability was done by, amongst many others, Hewitt and Savage [12], Heath and Sudderth [10], Diaconis and Freedman [8] and, in the context of the behavioural theory of imprecise probabilities that we are going to consider here, by Walley [19]. We refer to Kallenberg [14, 15] for modern, measure-theoretic discussions of exchangeability.

One of the reasons why exchangeability is deemed important, especially by Bayesians, is that, by virtue of de Finetti's representation theorem, an exchangeable model can be seen as a convex mixture of multinomial models. This has lent some support [2, 5, 7] to the claim that aleatory probabilities and i.i.d. processes can be eliminated from statistics and that we can restrict ourselves to exchangeable sequences instead; see Walley [19], Section 9.5.6 for a critical discussion of this claim.

De Finetti presented his study of exchangeability in terms of the behavioural notion of previsions, or fair prices. The central assumption underlying his approach is that a subject should be able to specify a fair price $P(f)$ for any risky transaction (which we will call a *gamble*) $f$ ([7], Chapter 3). This may not always be realistic, so it has been suggested that we should explicitly allow for a subject's indecision, by distinguishing between his *lower prevision* $\underline{P}(f)$, which is the supremum price for which he is willing to buy the gamble $f$, and his *upper prevision* $\overline{P}(f)$, which is the infimum price for which he is willing to sell $f$. For any real number $r$ strictly between $\underline{P}(f)$ and $\overline{P}(f)$, the subject is then not specifying a choice between selling or buying the gamble $f$ for $r$. Such lower and upper previsions are also subject to certain rationality or coherence criteria, in very much the same way that (precise) previsions are, in de Finetti's account. The resulting *theory of coherent lower previsions*, brilliantly defended by Walley [19], generalises de Finetti's behavioural treatment of subjective, epistemic probability and is briefly overviewed in Section 2.

Also, in this theory, it is interesting to consider the consequences of a subject's *exchangeability assessment*, that is, that the order in which we consider a number of random variables has no impact. This is our motivation for studying exchangeable *lower* previsions in this paper. An assessment of exchangeability will have a clear impact on the structure of so-called *exchangeable* coherent lower previsions. We will show that such a prevision can be written as a combination of (i) a coherent (linear) prevision expressing that permutations of realisations of such sequences are considered equally likely, and (ii) a coherent lower prevision for the 'frequency' of occurrence of the different values the random variables can take. Of course, this is the essence of representation in de Finetti's sense – we generalise his results to coherent lower previsions.

Before we go on, we want to draw attention to a number of distinctive features of our approach. First, the usual proofs of the representation theorem, such as the ones given by de Finetti [5], Heath and Sudderth [10] and Kallenberg [15], do not lend themselves very easily to generalisation in terms of coherent lower previsions. In principle, it would be possible, at least in some cases, to start with the versions already known for (precise) previsions and to derive their counterparts for lower previsions using so-called lower envelope theorems; see Section 2 for more details. This is the method that Walley [19], Sections 9.5.3 and 9.5.4, suggests. However, we have decided to follow a different route: we derive our results directly for lower previsions, using an approach based on Bernstein polynomials, and we obtain the ones for previsions as special cases. We believe this method to be more elegant and self-contained, and it certainly has the additional benefit of drawing attention to what we feel is the essence of de Finetti's representation theorem: specifying a coherent belief model for a countable exchangeable sequence is tantamount to specifying a coherent (lower) prevision on the linear space of polynomials on some simplex, and nothing more.

Second, we will focus on – and use the language of – (lower and upper) previsions for gambles, rather than (lower and upper) probabilities for events: in the behavioural theory of imprecise probabilities, the language of gambles is much more expressive than that of events and we need its full expressive power to derive our results.

The paper is organised as follows. In Section 2, we introduce a number of results from the theory of coherent lower previsions necessary to understand the rest of the paper. In Section 3, we define exchangeability for finite sequences of random variables and establish a representation of coherent exchangeable lower previsions in terms of sampling without replacement. In Section 4, we extend the notion of exchangeability to countable sequences of random variables and in Section 5 we generalise de Finetti's representation theorem (in terms of multinomial sampling) to exchangeable coherent lower previsions. In the Appendix, we have gathered a few useful results about Bernstein polynomials.

## 2. Lower previsions, random variables and their distributions

In this section, we provide a brief summary of ideas and results from the theory of coherent lower previsions [19].

### 2.1. Epistemic uncertainty models

Consider a *random variable* $X$ that may assume values $x$ in some non-empty set $\mathcal{X}$. By 'random', we mean that a subject is uncertain about the actual value of the variable $X$, that is, does not know what this actual value is.

Our subject may entertain certain beliefs about the value of $X$. We try and model his beliefs mathematically using the concept of a *gamble* on $\mathcal{X}$, which is a bounded map $f$ from $\mathcal{X}$ to the set $\mathbb{R}$ of real numbers. We denote by $\mathcal{L}(\mathcal{X})$ the set of all gambles on $\mathcal{X}$.

De Finetti [7] proposed the modelling of a subject's beliefs by eliciting his fair price, or *prevision*, $P(f)$ for certain gambles $f$. This $P(f)$ can be defined as the unique real number $p$ such that the subject is willing to buy the gamble $f$ for all prices $s$ (that is, accept the gamble $f - s$) and sell $f$ for all prices $t$ (that is, accept the gamble $t - f$) for all $s < p < t$. The problem with this approach is that it presupposes that there is such a real number, or, in other words, that the subject, whatever his beliefs about $X$ are, is willing, for (almost) every real $r$, to make a choice between buying $f$ for the price $r$ or selling it for that price.

### 2.2. Coherent lower previsions and natural extension

A way to address this problem is to consider a model that allows our subject to be undecided for some prices $r$. This is done in Walley's [19] theory of lower and upper previsions. The *lower prevision* of the gamble $f$, $\underline{P}(f)$, is our subject's supremum acceptable buying price for $f$; similarly, our subject's *upper prevision*, $\overline{P}(f)$, is his infimum acceptable selling price for $f$. Hence, he is willing to buy the gamble $f$ for all prices $s < \underline{P}(f)$ and sell $f$ for all prices $t > \overline{P}(f)$, but he may be undecided for prices $\underline{P}(f) \leq p \leq \overline{P}(f)$.

Since buying the gamble $f$ for a price $s$ is the same as selling the gamble $-f$ for the price $-s$, the lower and upper previsions are *conjugate* functions: $\underline{P}(f) = -\overline{P}(-f)$ for any gamble $f$. This allows us to concentrate on one of them since we can immediately derive results for the other. In this paper, we focus mainly on lower previsions.

The *lower probability* $\underline{P}(A)$ of an event $A \subseteq \mathcal{X}$ is defined as the lower prevision of its indicator $I_A$: $\underline{P}(A) = \underline{P}(I_A)$; $I_A$ is the gamble that assumes the value one on $A$ and zero elsewhere. For the *upper probability* $\overline{P}(A)$ of $A$, we similarly have that $\overline{P}(A) = \overline{P}(I_A)$.

For lower previsions, the most important rationality criterion is that of coherence. If a lower prevision $\underline{P}$ is defined on a linear space of gambles $\mathcal{K}$, then it turns out to be *coherent* if and only if it satisfies the following conditions. For any gambles $f$ and $g$ in $\mathcal{K}$ and any non-negative real number $\lambda$, it should hold that:

(P1) $\underline{P}(f) \geq \inf f$ [accepting sure gains];
(P2) $\underline{P}(\lambda f) = \lambda \underline{P}(f)$ [non-negative homogeneity];
(P3) $\underline{P}(f + g) \geq \underline{P}(f) + \underline{P}(g)$ [superadditivity].

The following special properties hold for a coherent lower prevision whenever the gambles involved belong to its domain:

(i) $\underline{P}$ is *monotone*, that is, if $f \leq g$, then $\underline{P}(f) \leq \underline{P}(g)$;
(ii) $\inf f \leq \underline{P}(f) \leq \overline{P}(f) \leq \sup f$.

Moreover, coherent lower and upper previsions are continuous with respect to uniform convergence of gambles.

### 2.3. Linear previsions

If the lower prevision $\underline{P}(f)$ and the upper prevision $\overline{P}(f)$ for a gamble $f$ happen to coincide, then the value $P(f) = \underline{P}(f) = \overline{P}(f)$ is called the subject's (precise) *prevision* for $f$. Previsions are fair prices in de Finetti's [7] sense. We shall call them *precise* probability models and lower previsions will be called *imprecise*.

A prevision on the set $\mathcal{L}(\mathcal{X})$ of all gambles is linear if and only if it is a positive ($f \geq 0 \Rightarrow P(f) \geq 0$) and normed ($P(1) = 1$) real linear functional. A prevision on a general domain is linear if and only if it can be extended to a linear prevision on all gambles. We shall denote by $\mathbb{P}(\mathcal{X})$ the set of all linear previsions on $\mathcal{L}(\mathcal{X})$.

There is an interesting link between precise and imprecise probability models, expressed via the so-called *lower envelope theorem* as follows. A lower prevision $\underline{P}$ on some domain $\mathcal{K}$ is coherent if and only if it is the *lower envelope* of some set of linear previsions and, in particular, of the convex set $\mathcal{M}(\underline{P})$ of all linear previsions that dominate it: for all $f$ in $\mathcal{K}$, $\underline{P}(f) = \inf\{P(f) : P \in \mathcal{M}(\underline{P})\}$, where $\mathcal{M}(\underline{P}) := \{P \in \mathbb{P}(\mathcal{X}) : (\forall f \in \mathcal{K})(P(f) \geq \underline{P}(f))\}$.

### 2.4. The distribution of a random variable

We call a subject's coherent lower prevision $\underline{P}$ on $\mathcal{L}(\mathcal{X})$, modelling his beliefs about the value that a random variable $X$ assumes in the set $\mathcal{X}$, his *distribution* for that random variable.

If we now consider another set $\mathcal{Y}$ and a map $\varphi$ from $\mathcal{X}$ to $\mathcal{Y}$, then we can consider $Y := \varphi(X)$ as a random variable assuming values in $\mathcal{Y}$. With a gamble $h$ on $\mathcal{Y}$, there corresponds a gamble $h \circ \varphi$ on $\mathcal{X}$ whose lower prevision is $\underline{P}(h \circ \varphi)$. This leads us to define the distribution of $Y = \varphi(X)$ as the *induced* coherent lower prevision $\underline{Q}$ on $\mathcal{L}(\mathcal{Y})$, defined by

$$\underline{Q}(h) := \underline{P}(h \circ \varphi), \qquad h \in \mathcal{L}(\mathcal{Y}).$$

This notion generalises that of an induced probability measure.

Finally, consider a sequence of random variables $X_n$, all taking values in some metric space $S$. Denote by $\mathcal{C}(S)$ the set of all continuous gambles on $S$. For each random variable $X_n$, we have a distribution in the form of a coherent lower prevision $\underline{P}_{X_n}$ on $\mathcal{L}(S)$. We then say that the random variables *converge in distribution* if for all $h \in \mathcal{C}(S)$, the sequence of real numbers $\underline{P}_{X_n}(h)$ converges to some real number, which we denote by $\underline{P}(h)$. The limit lower prevision $\underline{P}$ on $\mathcal{C}(S)$ that we can define in this way is coherent, because a pointwise limit of coherent lower previsions always is.

## 3. Exchangeable random variables

We are now ready to recall Walley's [19], Section 9.5, notion of exchangeability in the context of the theory of coherent lower previsions. We shall see that it generalises de Finetti's definition for linear previsions [5, 7].

### 3.1. Definition and basic properties

Consider $N \geq 1$ random variables $X_1, \ldots, X_N$ taking values in a non-empty and finite set $\mathcal{X}$. A subject's beliefs about the values that these random variables $\mathbf{X} = (X_1, \ldots, X_N)$ assume jointly in $\mathcal{X}^N$ is given by their (joint) distribution, which is a coherent lower prevision $\underline{P}^N$ defined on the set $\mathcal{L}(\mathcal{X}^N)$.

Let us denote by $\mathcal{P}_N$ the set of all permutations of $\{1, \ldots, N\}$. With any such permutation $\pi$, we can associate, by the procedure of lifting, a permutation of $\mathcal{X}^N$, also denoted by $\pi$, that maps any $\mathbf{x} = (x_1, \ldots, x_N)$ in $\mathcal{X}^N$ to $\pi \mathbf{x} := (x_{\pi(1)}, \ldots, x_{\pi(N)})$. Similarly, with any gamble $f$ on $\mathcal{X}^N$, we can consider the permuted gamble $\pi f := f \circ \pi$.

A subject judges the random variables $X_1, \ldots, X_N$ to be *exchangeable* when he is disposed to exchange any gamble $f$ for the permuted gamble $\pi f$, meaning that $\underline{P}^N(\pi f - f) \geq 0$, for any permutation $\pi$. Taking into account the properties of coherence, this means that

$$\overline{P}^N(f - \pi f) = \underline{P}^N(f - \pi f) = \overline{P}^N(\pi f - f) = \underline{P}^N(\pi f - f) = 0$$

for all gambles $f$ on $\mathcal{X}^N$ and all permutations $\pi$ in $\mathcal{P}_N$. In this case, we also call the joint coherent lower prevision $\underline{P}^N$ *exchangeable*. A subject will make an assumption of exchangeability when there is evidence that the processes generating the values of the random variables are (physically) similar [19], Section 9.5.2, and consequently the order in which the variables are observed is not important.

When $\underline{P}^N$ is, in particular, a linear prevision $P^N$, exchangeability is equivalent to having $P^N(\pi f) = P^N(f)$ for all gambles $f$ and all permutations $\pi$. The following proposition, mentioned by Walley [19], Section 9.5, and whose proof is immediate and therefore omitted, establishes an even stronger link between Walley's and de Finetti's notions of exchangeability.

**Proposition 1.** *A coherent lower prevision $\underline{P}^N$ is exchangeable if and only if all the linear previsions $P^N$ in $\mathcal{M}(\underline{P}^N)$ are exchangeable.*

Clearly, if $X_1, \ldots, X_N$ are exchangeable, then any permutation $X_{\pi(1)}, \ldots, X_{\pi(N)}$ is also exchangeable and has the same distribution $\underline{P}^N$. Moreover, any selection of $1 \leq n \leq N$ random variables from amongst the $X_1, \ldots, X_N$ are exchangeable too and their distribution is given by $\underline{P}^n$, which is the $\mathcal{X}^n$-marginal of $\underline{P}^N$, defined by $\underline{P}^n(f) := \underline{P}^N(\widetilde{f})$ for all gambles $f$ on $\mathcal{X}^n$, where the gamble $\widetilde{f}$ on $\mathcal{X}^N$ is the *cylindrical extension* of $f$ to $\mathcal{X}^N$, given by $\widetilde{f}(z_1, \ldots, z_N) := f(z_1, \ldots, z_n)$ for all $(z_1, \ldots, z_N)$ in $\mathcal{X}^N$.

### 3.2. Count vectors

Interestingly, exchangeable coherent lower previsions have a very simple representation, in terms of sampling without replacement. To see how this comes about, consider any $\mathbf{x} \in \mathcal{X}^N$. The so-called (permutation) *invariant atom*

$$[\mathbf{x}] := \{\pi \mathbf{x} : \pi \in \mathcal{P}_N\}$$

is then the smallest non-empty subset of $\mathcal{X}^N$ that contains $\mathbf{x}$ and is invariant under all permutations $\pi$ in $\mathcal{P}_N$. We shall denote the set of permutation invariant atoms of $\mathcal{X}^N$ by $\mathcal{A}^N$. This constitutes a partition of the set $\mathcal{X}^N$. We can characterise these invariant atoms using the *counting maps* $T_x^N : \mathcal{X}^N \to \mathbb{N}_0$ defined for all $x$ in $\mathcal{X}$ in such a way that

$$T_x^N(\mathbf{z}) = T_x^N(z_1, \ldots, z_N) := |\{k \in \{1, \ldots, N\} : z_k = x\}|$$

is the number of components of the $N$-tuple $\mathbf{z}$ that assume the value $x$. Here, $|A|$ denotes the number of elements in a finite set $A$ and $\mathbb{N}_0$ is the set of all non-negative integers (including zero). We shall denote by $\mathbf{T}^N$ the vector-valued map from $\mathcal{X}^N$ to $\mathbb{N}_0^{\mathcal{X}}$ whose component maps are the $T_x^N$, $x \in \mathcal{X}$. $\mathbf{T}^N$ actually assumes values in the set of *count vectors*

$$\mathcal{N}^N := \left\{\mathbf{m} \in \mathbb{N}_0^{\mathcal{X}} : \sum_{x \in \mathcal{X}} m_x = N\right\}.$$

The counting map $\mathbf{T}^N$ can be interpreted as a bijection (one-to-one and onto) between the set of invariant atoms $\mathcal{A}^N$ and the set of count vectors $\mathcal{N}^N$, and we can identify any invariant atom $[\mathbf{z}]$ by the count vector $\mathbf{m} = \mathbf{T}^N(\mathbf{z})$ of any (and therefore all) of its elements. We therefore also denote this atom by $[\mathbf{m}]$. Clearly $\mathbf{y} \in [\mathbf{m}]$ if and only if $\mathbf{T}^N(\mathbf{y}) = \mathbf{m}$. The number of elements $\nu(\mathbf{m})$ in any invariant atom $[\mathbf{m}]$ is given by

$$\nu(\mathbf{m}) := \binom{N}{\mathbf{m}} = \frac{N!}{\prod_{x \in \mathcal{X}} m_x!}.$$

If the joint random variable $\mathbf{X} = (X_1, \ldots, X_N)$ assumes the value $\mathbf{z}$ in $\mathcal{X}^N$, then the corresponding count vector assumes the value $\mathbf{T}^N(\mathbf{z})$ in $\mathcal{N}^N$. This means that we can see $\mathbf{T}^N(\mathbf{X}) = \mathbf{T}^N(X_1, \ldots, X_N)$ as a random variable in $\mathcal{N}^N$. If the available information about the values that $\mathbf{X}$ assumes in $\mathcal{X}^N$ is given by the coherent exchangeable lower prevision $\underline{P}^N$ (the distribution of $\mathbf{X}$), then the corresponding uncertainty model for the values that $\mathbf{T}^N(\mathbf{X})$ assumes in $\mathcal{N}^N$ is given by the coherent *induced* lower prevision $\underline{Q}^N$ on $\mathcal{L}(\mathcal{N}^N)$ (the distribution of $\mathbf{T}^N(\mathbf{X})$), given by

$$\underline{Q}^N(h) := \underline{P}^N(h \circ \mathbf{T}^N) = \underline{P}^N \left( \sum_{\mathbf{m} \in \mathcal{N}^N} h(\mathbf{m}) I_{[\mathbf{m}]} \right) \qquad \text{for all gambles } h \text{ on } \mathcal{N}^N. \quad (1)$$

We now come to a theorem showing that, conversely, any exchangeable coherent lower prevision $\underline{P}^N$ is in fact *completely determined* by the corresponding distribution $\underline{Q}^N$ of the count vectors, also called its *count distribution*.

Consider an urn with $N$ balls of different types, where the different types are characterised by the elements $x$ of the set $\mathcal{X}$. Suppose the *composition* of the urn is given by the count vector $\mathbf{m} \in \mathcal{N}^N$, meaning that $m_x$ balls are of type $x$ for $x \in \mathcal{X}$. We are now going to subsequently select (in a random way) $N$ balls from the urn, without replacing them. It follows that for any gamble $f$ on $\mathcal{X}^N$, its (precise) prevision (or expectation) is given by

$$MuHy^N(f|\mathbf{m}) := \frac{1}{\nu(\mathbf{m})} \sum_{\mathbf{z} \in [\mathbf{m}]} f(\mathbf{z}).$$

The linear prevision $MuHy^N(\cdot|\mathbf{m})$ is the one associated with a *multiple hypergeometric distribution* ([13], Chapter 39), whence the notation. For any permutation $\pi$ of $\{1, \ldots, N\}$,

$$MuHy^N(\pi f|\mathbf{m}) = \frac{1}{\nu(\mathbf{m})} \sum_{\mathbf{z} \in [\mathbf{m}]} f(\pi \mathbf{z}) = \frac{1}{\nu(\mathbf{m})} \sum_{\pi^{-1}\mathbf{z} \in [\mathbf{m}]} f(\mathbf{z}) = MuHy^N(f|\mathbf{m})$$

since $\pi^{-1}\mathbf{z} \in [\mathbf{m}]$ if and only if $\mathbf{z} \in [\mathbf{m}]$. This means that the linear prevision $MuHy^N(\cdot|\mathbf{m})$ is exchangeable. The following theorem establishes an even stronger result. It is an immediate consequence of a much more general representation result by de Cooman and Miranda [4], Theorem 30.

**Theorem 2 (Representation theorem for finite sequences of exchangeable variables).** *Let $N \geq 1$. A coherent lower prevision $\underline{P}^N$ on $\mathcal{L}(\mathcal{X}^N)$ is exchangeable if and only if it there is some coherent lower prevision $\underline{Q}$ on $\mathcal{L}(\mathcal{N}^N)$ such that $\underline{P}^N(f) =$* $\underline{Q}(MuHy^N(f|\cdot))$ *for all gambles $f$ on $\mathcal{X}^N$. If a coherent lower prevision $\underline{P}^N$ on $\mathcal{L}(\mathcal{X}^N)$ is exchangeable, then the corresponding $\underline{Q}$ is given by equation (1).*

This theorem implies that any collection of $N$ exchangeable random variables in $\mathcal{X}$ can be seen as the result of $N$ random draws without replacement from an urn with $N$ balls whose types are characterised by the elements $x$ of $\mathcal{X}$ and whose composition **m** is unknown, but for which the available information about the composition is modelled by a coherent lower prevision on $\mathcal{L}(\mathcal{N}^N)$.[1]

That exchangeable linear previsions can be interpreted in terms of sampling without replacement from an urn with unknown composition is of course well known and essentially goes back to de Finetti's work on exchangeability [1, 5]. Heath and Sudderth [10] give a simple proof for variables that may assume two values. However, we believe our proof of the much more general representation result ([4], Theorem 30), to be conceptually even simpler than Heath and Sudderth's proof.

## 4. Exchangeable sequences

### 4.1. Definitions

Consider a countable sequence $X_1, \ldots, X_n, \ldots$ of random variables taking values in the same non-empty set $\mathcal{X}$. This sequence is called *exchangeable* if any finite collection of random variables taken from this sequence is exchangeable.

We can also consider the exchangeable sequence as a single random variable **X** assuming values in the set $\mathcal{X}^\mathbb{N}$, where $\mathbb{N}$ is the set of natural numbers (positive integers, without zero). Its possible values **x** are sequences $x_1, \ldots, x_n, \ldots$ of elements of $\mathcal{X}$ or, in other words, maps from $\mathbb{N}$ to $\mathcal{X}$. We can model the available information about the value that **X** assumes in $\mathcal{X}^\mathbb{N}$ by a coherent lower prevision $\underline{P}^\mathbb{N}$ on $\mathcal{L}(\mathcal{X}^\mathbb{N})$, called the *distribution* of the exchangeable random sequence **X**.

The random sequence **X**, or its distribution $\underline{P}^\mathbb{N}$, is clearly exchangeable if and only if all of its $\mathcal{X}^n$-*marginals* $\underline{P}^n$ are exchangeable for $n \geq 1$. These marginals $\underline{P}^n$ on $\mathcal{L}(\mathcal{X}^n)$ are defined as follows: for any gamble $f$ on $\mathcal{X}^n$, $\underline{P}^n(f) := \underline{P}^\mathbb{N}(\widetilde{f})$, where $\widetilde{f}$ is the cylindrical extension of $f$ to $\mathcal{X}^\mathbb{N}$ defined by $\widetilde{f}(\mathbf{x}) := f(x_1, \ldots, x_n)$ for all $\mathbf{x} = (x_1, \ldots, x_n, x_{n+1}, \ldots)$ in $\mathcal{X}^\mathbb{N}$. In addition, the family of exchangeable coherent lower previsions $\underline{P}^n$, $n \geq 1$, satisfies the '*time consistency*' requirement

$$\underline{P}^n(f) = \underline{P}^{n+k}(\widetilde{f}) \tag{2}$$

for all $n \geq 1$, $k \geq 0$ and all gambles $f$ on $\mathcal{X}^n$, where $\widetilde{f}$ now denotes the cylindrical extension of $f$ to $\mathcal{X}^{n+k}$: $\underline{P}^n$ should be the $\mathcal{X}^n$-marginal of any $\underline{P}^{n+k}$.

It follows at once that any finite collection of $n \geq 1$ random variables taken from such an exchangeable sequence has the same distribution as the first $n$ variables $X_1, \ldots, X_n$, which is the exchangeable coherent lower prevision $\underline{P}^n$ on $\mathcal{L}(\mathcal{X}^n)$.

---

[1]Walley [19], Chapter 9, also mentions this result for exchangeable coherent lower previsions.

Conversely, suppose we have a collection of exchangeable coherent lower previsions $\underline{P}^n$ on $\mathcal{L}(\mathcal{X}^n)$, $n \geq 1$, that satisfy the time consistency requirement (2). Then any coherent lower prevision $\underline{P}^{\mathbb{N}}$ on $\mathcal{L}(\mathcal{X}^{\mathbb{N}})$ that has $\mathcal{X}^n$-marginals $\underline{P}^n$ is exchangeable. The smallest, or most conservative, such (exchangeable) coherent lower prevision is given by

$$\underline{E}^{\mathbb{N}}(f) := \sup_{n \in \mathbb{N}} \underline{P}^n(\underline{\text{proj}}_n(f)) = \lim_{n \to \infty} \underline{P}^n(\underline{\text{proj}}_n(f)),$$

where $f$ is any gamble on $\mathcal{X}^{\mathbb{N}}$ and its *lower projection* $\underline{\text{proj}}_n(f)$ on $\mathcal{X}^n$ is the gamble on $\mathcal{X}^n$ that is defined by $\underline{\text{proj}}_n(f)(\mathbf{x}) := \inf_{\mathbf{z} \in \mathcal{X}^{\mathbb{N}} : z_k = x_k, k=1,\ldots,n} f(\mathbf{z})$ for all $\mathbf{x} \in \mathcal{X}^n$; see de Cooman and Miranda [3], Section 5, for more details.

## 4.2. Time consistency of the count distributions

It is of crucial interest for what follows to determine the consequences of the time consistency requirement (2) on the marginals $\underline{P}^n$ for the corresponding family $\underline{Q}^n$, $n \geq 1$, of distributions of the count vectors $\mathbf{T}^n(X_1, \ldots, X_n)$. Consider, therefore, $n \geq 1$, $k \geq 0$ and any gamble $h$ on $\mathcal{N}^n$. If we let $f := h \circ \mathbf{T}^n$, then

$$\underline{Q}^n(h) = \underline{P}^n(f) = \underline{P}^{n+k}(\widetilde{f}) = \underline{Q}^{n+k}(MuHy^{n+k}(\widetilde{f}|\cdot)),$$

where the first equality follows from equation (1), the second from equation (2) and the last from Theorem 2. Now, for any $\mathbf{m}'$ in $\mathcal{N}^{n+k}$ and any $\mathbf{z}' = (\mathbf{z}, \mathbf{y})$ in $\mathcal{X}^{n+k} = \mathcal{X}^n \times \mathcal{X}^k$, we have that $\mathbf{T}^{n+k}(\mathbf{z}') = \mathbf{T}^n(\mathbf{z}) + \mathbf{T}^k(\mathbf{y})$ and therefore

$$MuHy^{n+k}(\widetilde{f}|\mathbf{m}') = \sum_{\mathbf{m} \in \mathcal{N}^n} \frac{\nu(\mathbf{m}' - \mathbf{m})\nu(\mathbf{m})}{\nu(\mathbf{m}')} h(\mathbf{m}),$$

taking into account that $MuHy^n(f|\mathbf{m}) = h(\mathbf{m})$ and that $\nu(\mathbf{m}' - \mathbf{m})$ is zero unless $\mathbf{m} \leq \mathbf{m}'$. So we see that time consistency is equivalent to

$$\underline{Q}^n(h) = \underline{Q}^{n+k}\left(\sum_{\mathbf{m} \in \mathcal{N}^n} \frac{\nu(\cdot - \mathbf{m})\nu(\mathbf{m})}{\nu(\cdot)} h(\mathbf{m})\right) \tag{3}$$

for all $n \geq 1$, $k \geq 0$ and $h \in \mathcal{L}(\mathcal{N}^n)$.

# 5. A representation theorem for exchangeable sequences

De Finetti [5, 7] has proven a representation result for exchangeable sequences with linear previsions that generalises Theorem 2 and where multinomial distributions take over the role that the multiple hypergeometric ones play for finite collections of exchangeable variables. One simple and intuitive way (see also [7], p. 218) to understand why the representation result can be thus extended from finite collections to countable sequences, is based on the fact that the multinomial distribution can be seen as a limit of multiple hypergeometric ones ([13], Chapter 39). This is also the central idea behind Heath and Sudderth [10] simple proof of this representation result in the case of variables that may only assume two possible values.

However, there is another, arguably even simpler, approach to proving the same results, which we present here. It also works for exchangeability in the context of coherent lower previsions. And, as we shall have occasion to explain further on, it has the additional advantage of clearly indicating what the 'representation' is and where it is uniquely defined.

### 5.1. Multinomial processes are exchangeable

Consider a sequence of i.i.d. random variables $Y_1, \ldots, Y_n, \ldots$ with common probability mass function $\boldsymbol{\theta}$: the probability that $Y_n = x$ is $\theta_x$ for $x \in \mathcal{X}$. Observe that $\boldsymbol{\theta}$ is an element of the $\mathcal{X}$-*simplex*

$$\Sigma = \left\{ \boldsymbol{\theta} \in \mathbb{R}^{\mathcal{X}} : (\forall x \in \mathcal{X}) \, (\theta_x \geq 0) \text{ and } \sum_{x \in \mathcal{X}} \theta_x = 1 \right\}.$$

Then, for any $n \geq 1$ and any $\mathbf{z}$ in $\mathcal{X}^n$, the probability that $(Y_1, \ldots, Y_n)$ is equal to $\mathbf{z}$ is given by $\prod_{x \in \mathcal{X}} \theta_x^{T_x(\mathbf{z})}$, which yields the *multinomial mass function* ([13], Chapter 35). As a result, we have for any gamble $f$ on $\mathcal{X}^n$ that its corresponding (multinomial) prevision (expectation) is given by

$$Mn^n(f|\boldsymbol{\theta}) = CoMn^n(MuHy^n(f|\cdot)|\boldsymbol{\theta}), \tag{4}$$

where we defined the (count multinomial) linear prevision $CoMn^n(\cdot|\boldsymbol{\theta})$ on $\mathcal{L}(\mathcal{N}^n)$ by

$$CoMn^n(g|\boldsymbol{\theta}) = \sum_{\mathbf{m} \in \mathcal{N}^n} g(\mathbf{m}) \nu(\mathbf{m}) \prod_{x \in \mathcal{X}} \theta_x^{m_x}, \tag{5}$$

where $g$ is any gamble on $\mathcal{N}^n$. The corresponding probability mass for any count vector $\mathbf{m}$,

$$CoMn^n(\{\mathbf{m}\}|\boldsymbol{\theta}) = \nu(\mathbf{m}) \prod_{x \in \mathcal{X}} \theta_x^{m_x} =: B_{\mathbf{m}}(\boldsymbol{\theta}) \tag{6}$$

is the probability of observing some value $\mathbf{z}$ for $(Y_1, \ldots, Y_n)$ whose count vector is $\mathbf{m}$. The polynomial function $B_{\mathbf{m}}$ on the $\mathcal{X}$-simplex is called a (multivariate) *Bernstein (basis) polynomial*. The set $\{B_{\mathbf{m}} : \mathbf{m} \in \mathcal{N}^n\}$ of all Bernstein (basis) polynomials of fixed degree $n$ forms a basis for the linear space of all (multivariate) polynomials on $\Sigma$ whose degree is at most $n$, hence their name. If we have a polynomial $p$ of degree $m$, this means that for any $n \geq m$, $p$ has a *unique* (Bernstein) decomposition $b_p^n \in \mathcal{L}(\mathcal{N}^n)$ such that $p = \sum_{\mathbf{m} \in \mathcal{N}^n} b_p^n(\mathbf{m}) B_{\mathbf{m}}$. If we combine this with equations (5) and (6), we find that $b_p^n$ is the unique gamble on $\mathcal{N}^n$ such that $CoMn^n(b_p^n | \cdot) = p$.

We deduce from equation (4) and Theorem 2 that the linear prevision $Mn^n(\cdot|\boldsymbol{\theta})$ on $\mathcal{L}(\mathcal{X}^n)$ is exchangeable and that $CoMn^n(\cdot|\boldsymbol{\theta})$ is the corresponding distribution for the corresponding count vectors $\mathbf{T}^n(Y_1,\ldots,Y_n)$. Therefore, the sequence of i.i.d. random variables $Y_1,\ldots,Y_n,\ldots$ is exchangeable.

### 5.2. A representation theorem

Consider the linear subspace of $\mathcal{L}(\Sigma)$,

$$\mathcal{V}(\Sigma) := \{CoMn^n(g|\cdot): n \geq 1, g \in \mathcal{L}(\mathcal{N}^n)\} = \{Mn^n(f|\cdot): n \geq 1, f \in \mathcal{L}(\mathcal{X}^n)\},$$

each of whose elements is a *polynomial function* on the $\mathcal{X}$-simplex,

$$CoMn^n(g|\boldsymbol{\theta}) = \sum_{\mathbf{m} \in \mathcal{N}^n} g(\mathbf{m})\nu(\mathbf{m}) \prod_{x \in \mathcal{X}} \theta_x^{m_x} = \sum_{\mathbf{m} \in \mathcal{N}^n} g(\mathbf{m}) B_{\mathbf{m}}(\boldsymbol{\theta}),$$

and is actually a linear combination of Bernstein basis polynomials $B_{\mathbf{m}}$ with coefficients $g(\mathbf{m})$. So, $\mathcal{V}(\Sigma)$ is the linear space spanned by all Bernstein basis polynomials and is therefore the set of all polynomials on the $\mathcal{X}$-simplex $\Sigma$.

Now, if $\underline{R}$ is any coherent lower prevision on $\mathcal{L}(\Sigma)$, then it is easy to see that the family of coherent lower previsions $\underline{P}^n$, $n \geq 1$, defined by

$$\underline{P}^n(f) = \underline{R}(Mn^n(f|\cdot)), \qquad f \in \mathcal{L}(\mathcal{X}^n), \tag{7}$$

is still exchangeable and time consistent, and that the corresponding count distributions are

$$\underline{Q}^n(g) = \underline{R}(CoMn^n(g|\cdot)), \qquad g \in \mathcal{L}(\mathcal{N}^n). \tag{8}$$

Here, we are going to show that a converse result also holds: for any time-consistent family of exchangeable coherent lower previsions $\underline{P}^n$, $n \geq 1$, there is a coherent lower prevision $\underline{R}$ on $\mathcal{V}(\Sigma)$ such that equation (7), or its reformulation for counts (8), holds. We call such an $\underline{R}$ a *representation*, or representing coherent lower prevision, for the family $\underline{P}^n$. Of course, any representing $\underline{R}$, if it exists, is uniquely determined on $\mathcal{V}(\Sigma)$.

So consider a family of coherent lower previsions $\underline{Q}^n$ on $\mathcal{L}(\mathcal{N}^n)$, $n \geq 1$, that are time consistent. It suffices to find an $\underline{R}$ such that (8) holds because the corresponding exchangeable lower previsions $\underline{P}^n$ on $\mathcal{L}(\mathcal{X}^n)$ are then uniquely determined by Theorem 2, and automatically satisfy the condition (7). Our proposal is to *define* the functional $\underline{R}$ on the set $\mathcal{V}(\Sigma)$ as follows: *consider any element $p$ of $\mathcal{V}(\Sigma)$. Then, by definition, there is some $n \geq 1$ and a corresponding unique $b_p^n \in \mathcal{L}(\mathcal{N}^n)$ such that $p = CoMn^n(b_p^n|\cdot)$. We then let $\underline{R}(p) := \underline{Q}^n(b_p^n)$.*

The first thing to check is whether this definition is consistent.

**Lemma 3.** *Let $p$ be a polynomial of degree $m$ and let $n_1, n_2 \geq m$. Then $\underline{Q}^{n_1}(b_p^{n_1}) = \underline{Q}^{n_2}(b_p^{n_2})$.*

**Proof.** We may assume without loss of generality that $n_2 \geq n_1$. The Bernstein decompositions $b_p^{n_1}$ and $b_p^{n_2}$ are then related by Zhou's formula [see equation (10) in the Appendix]:

$$b_p^{n_2}(\mathbf{m}_2) = \sum_{\mathbf{m}_1 \in \mathcal{N}^{n_1}} \frac{\nu(\mathbf{m}_2 - \mathbf{m}_1)\nu(\mathbf{m}_1)}{\nu(\mathbf{m}_2)} b_p^{n_1}(\mathbf{m}_1), \qquad \mathbf{m}_2 \in \mathcal{N}^{n_2}.$$

Consequently, by the time consistency requirement (3), $\underline{Q}^{n_2}(b_p^{n_2}) = \underline{Q}^{n_1}(b_p^{n_1})$. □

**Lemma 4.** $\underline{R}$ *is a coherent lower prevision on the linear space* $\mathcal{V}(\Sigma)$.

**Proof.** We show that $\underline{R}$ satisfies the necessary and sufficient conditions (P1)–(P3) for coherence of a lower prevision on a linear space.

We first prove that (P1) is satisfied. Consider any $p \in \mathcal{V}(\Sigma)$. Let $m$ be the degree of $p$. We must show that $\underline{R}(p) \geq \min p$. We find that $\underline{R}(p) = \underline{Q}^n(b_p^n) \geq \min b_p^n$ for all $n \geq m$ because of the coherence of $\underline{Q}^n$. However, equation (11) in the Appendix tells us that $\min b_p^n \uparrow \min p$, so we indeed have $\underline{R}(p) \geq \min p$.

Next, consider any $p$ in $\mathcal{V}(\Sigma)$ and any real $\lambda \geq 0$. Consider any $n$ that is not smaller than the degree of $p$. Since it is obvious that $b_{\lambda p}^n = \lambda b_p^n$, we get

$$\underline{R}(\lambda p) = \underline{Q}^n(b_{\lambda p}^n) = \underline{Q}^n(\lambda b_p^n) = \lambda \underline{Q}^n(b_p^n) = \lambda \underline{R}(p),$$

where the third equality follows from the coherence (non-negative homogeneity) of the count lower prevision $\underline{Q}^n$. This tells us that $\underline{R}$ satisfies (P2).

Finally, consider $p$ and $q$ in $\mathcal{V}(\Sigma)$, and any $n$ that is not smaller than the maximum of the degrees of $p$ and $q$. Since it is obvious that $b_{p+q}^n = b_p^n + b_q^n$, we get

$$\underline{R}(p + q) = \underline{Q}^n(b_{p+q}^n) = \underline{Q}^n(b_p^n + b_q^n) \geq \underline{Q}^n(b_p^n) + \underline{Q}^n(b_q^n) = \underline{R}(p) + \underline{R}(q),$$

where the inequality follows from the superadditivity of $\underline{Q}^n$. This tells us that $\underline{R}$ also satisfies (P3) and, as a consequence, it is coherent. □

We can summarise the argument above as follows.

**Theorem 5 (Representation theorem for exchangeable sequences).** *Given a time-consistent family of exchangeable coherent lower previsions $\underline{P}^n$ on $\mathcal{L}(\mathcal{X}^n)$, $n \geq 1$, there is a unique coherent lower prevision $\underline{R}$ on the linear space $\mathcal{V}(\Sigma)$ of all polynomial gambles on the $\mathcal{X}$-simplex such that for all $n \geq 1$, all $f \in \mathcal{L}(\mathcal{X}^n)$ and all $g \in \mathcal{L}(\mathcal{N}^n)$,*

$$\underline{P}^n(f) = \underline{R}(Mn^n(f|\cdot)) \quad and \quad \underline{Q}^n(g) = \underline{R}(CoMn^n(g|\cdot)). \tag{9}$$

Hence, the belief model governing any countable exchangeable sequence in $\mathcal{X}$ can be completely characterised by a coherent lower prevision on the linear space of polynomial gambles on $\Sigma$.

In the particular case where we have a time-consistent family of exchangeable *linear* previsions $P^n$ on $\mathcal{L}(\mathcal{X}^n), n \geq 1$, $\underline{R}$ will be a linear prevision $R$ on the linear space $\mathcal{V}(\Sigma)$ of all polynomial gambles on the $\mathcal{X}$-simplex. As such, it will be characterised by its values $\underline{R}(B_{\mathbf{m}})$ on the Bernstein basis polynomials $B_{\mathbf{m}}$, $\mathbf{m} \in \mathcal{N}^n$, $n \geq 1$, or on any other basis of $\mathcal{V}(\Sigma)$.

It is a consequence of coherence that $\underline{R}$ is also uniquely determined on the set $\mathcal{C}(\Sigma)$ of all continuous gambles on the $\mathcal{X}$-simplex $\Sigma$: by the Stone–Weierstrass theorem, any such gamble is the uniform limit of some sequence of polynomial gambles and coherence implies that the lower prevision of a uniform limit is the limit of the lower previsions.

This unicity result cannot be extended to more general (discontinuous) types of gambles: the coherent lower prevision $\underline{R}$ is not uniquely determined on the set of all gambles $\mathcal{L}(\Sigma)$ on the simplex and there may be different coherent lower previsions $\underline{R}^1$ and $\underline{R}^2$ on $\mathcal{L}(\Sigma)$ satisfying equation (9). But any such lower previsions will agree on the class $\mathcal{V}(\Sigma)$ of polynomial gambles, which is the class of gambles we need in order to characterise the exchangeable sequence.

We now investigate the meaning of the representing lower prevision $\underline{R}$ a bit further. Consider the sequence of so-called *frequency* random variables $\mathbf{F}_n := \mathbf{T}^n(X_1, \ldots, X_n)/n$ corresponding to an exchangeable sequence of random variables $X_1, \ldots, X_n, \ldots$ and assuming values in the $\mathcal{X}$-simplex $\Sigma$. The distribution $\underline{P}_{\mathbf{F}_n}$ of $\mathbf{F}_n$ is given by

$$\underline{P}_{\mathbf{F}_n}(h) := \underline{Q}^n\left(h \circ \frac{1}{n}\right) = \underline{R}\left(CoMn^n\left(h \circ \frac{1}{n}|\cdot\right)\right), \qquad h \in \mathcal{L}(\Sigma),$$

because we know that $\underline{Q}^n$ is the distribution of $\mathbf{T}^n(X_1, \ldots, X_n)$, and also taking into account Theorem 5 for the last equality. Now,

$$CoMn^n\left(h \circ \frac{1}{n}|\boldsymbol{\theta}\right) = \sum_{\mathbf{m} \in \mathcal{N}^n} h\left(\frac{\mathbf{m}}{n}\right) B_{\mathbf{m}}(\boldsymbol{\theta})$$

is the *Bernstein approximant* or *approximating Bernstein polynomial* of degree $n$ for the gamble $h$ and it is a known result (see [9], Section VII.2, or [11], Section 2) that the sequence of approximating Bernstein polynomials $CoMn^n(h \circ \frac{1}{n}|\cdot)$ converges uniformly to $h$ as $n \to \infty$ if $h$ is continuous. So, because $\underline{R}$ is uniquely defined and uniformly continuous on the set $\mathcal{C}(\Sigma)$, we find the following result.

**Theorem 6.** *For all continuous gambles $h$ on $\Sigma$, we have that*

$$\lim_{n \to \infty} \underline{P}_{\mathbf{F}_n}(h) = \underline{R}(h)$$

*or, in other words, the sequence of distributions $\underline{P}_{\mathbf{F}_n}$ converges pointwise to $\underline{R}$ on $\mathcal{C}(\Sigma)$ and, in this specific sense, the sample frequencies $\mathbf{F}_n$ converge in distribution.*

## 6. Conclusions

We have shown that the notion of exchangeability has a natural place in the theory of coherent lower previsions. Indeed, with our distinctive approach using Bernstein polyno-mials, and gambles rather than events, it seems fairly natural and easy to derive representation theorems directly for coherent lower previsions and to derive the corresponding results for precise probabilities (linear previsions) as special cases.

Interesting results can also be obtained in a context of predictive inference, where a coherent exchangeable lower prevision for $n+k$ variables is updated with the information that the first $n$ variables have been observed to assume certain values. For a fairly detailed discussion of these issues, we refer to de Cooman and Miranda [4], Section 9.3.

## Appendix: Multivariate Bernstein polynomials

To any $n \geq 0$ and $\mathbf{m} \in \mathcal{N}^n$, there corresponds a Bernstein (basis) polynomial of degree $n$ on $\Sigma$, given by $B_\mathbf{m}(\boldsymbol{\theta}) = \nu(\mathbf{m}) \prod_{x \in \mathcal{X}} \theta_x^{m_x}$, $\boldsymbol{\theta} \in \Sigma$. These polynomials have a number of very interesting properties (see, for instance, [17], Chapters 10 and 11):

- (B1) they are non-negative, and strictly positive in the interior of $\Sigma$;
- (B2) the set $\{B_\mathbf{m} : \mathbf{m} \in \mathcal{N}^n\}$ of all Bernstein polynomials of fixed degree $n$ forms a basis for the linear space of all polynomials whose degree is at most $n$.

Hence, for any polynomial $p$ of degree $m$, there is a unique gamble $b_p^n$ on $\mathcal{N}^n$ such that

$$p = \sum_{\mathbf{m} \in \mathcal{N}^n} b_p^n(\mathbf{m}) B_\mathbf{m} = CoMn^n(b_p^n | \cdot).$$

This tells us that each $p(\boldsymbol{\theta})$ is a convex combination of the Bernstein coefficients $b_p^n(\mathbf{m})$, $\mathbf{m} \in \mathcal{N}^n$, so $\min b_p^n \leq \min p \leq p(\boldsymbol{\theta}) \leq \max p \leq \max b_p^n$. It also follows that for all $k \geq 0$ and all $\boldsymbol{\mu}$ in $\mathcal{N}^{n+k}$,

$$b_p^{n+k}(\boldsymbol{\mu}) = \sum_{\mathbf{m} \in \mathcal{N}^n} \frac{\nu(\mathbf{m})\nu(\boldsymbol{\mu} - \mathbf{m})}{\nu(\boldsymbol{\mu})} b_p^n(\mathbf{m}). \tag{10}$$

This is *Zhou's formula* (see [17], Section 11.9). Moreover, since for any polynomial $p$ on $\Sigma$ of degree $m$, the $b_p^n$ converge uniformly to $p$ as $n \to \infty$ (see, for instance, [18]), it follows that

$$\lim_{\substack{n \to \infty \\ n \geq m}} [\min b_p^n, \max b_p^n] = [\min p, \max p] = p(\Sigma). \tag{11}$$

## Acknowledgements

We wish to acknowledge financial support in the form of research Grant G.0139.01 of the Flemish Fund for Scientific Research (FWO) and project TIN2008-06796-C04-01. The second author's research was financed by a Ph.D. grant from the Institute for the Promotion of Innovation through Science and Technology in Flanders (IWT Vlaanderen).

We would like to thank Jürgen Garloff for very helpful comments and pointers to the literature concerning multivariate Bernstein polynomials.